\documentclass[12pt]{amsart}
\usepackage{amsmath}
\usepackage{amsfonts,amssymb,amsthm}
\usepackage{graphicx}

\newcommand{\beq}{\begin{equation}}
\newcommand{\eeq}{\end{equation}}
\newcommand{\bbar}{\begin{eqnarray}}
\newcommand{\eear}{\end{eqnarray}}

\newcommand{\thm}[2]{\begin{#1} #2 \end{#1}}

\newtheorem{theorem}{Theorem}[section]
\newtheorem{itheorem}{Theorem}[section]
\newtheorem{lemma}[theorem]{Lemma}
\newtheorem{ilemma}[itheorem]{Lemma}
\newtheorem{corollary}[theorem]{Corollary}

\newtheorem{remark}[theorem]{Remark}

\newtheorem{example}[theorem]{Example}
\newtheorem{definition}[theorem]{Definition}
\newtheorem{idefinition}[itheorem]{Definition}

\begin{document}

\title{The moment zeta function and applications}

\author{Igor Rivin}
\address{Mathematics department, University of Manchester,
Oxford Road, Manchester M13 9PL, UK}
\address{Mathematics Department, Temple University,
Philadelphia, PA 19122}
\address{Mathematics Department, Princeton University, Princeton,
NJ 08544}
\email{irivin@math.princeton.edu} \thanks{The author would like 
to think the EPSRC and the NSF for support, and Natalia Komarova 
and  Ilan Vardi for useful conversations. }

\begin{abstract}
Motivated by a probabilistic analysis of a simple game (itself
inspired by a problem in computational learning theory) we introduce
the \emph{moment zeta function} of a probability distribution, and
study in depth some asymptotic properties of the moment zeta function
of those distributions supported in the interval $[0, 1].$ One example
of such zeta functions is Riemann's zeta function (which is the moment
zeta function of the uniform distribution in $[0, 1].$ For Riemann's
zeta function we are able to show particularly sharp versions of our
results.
\end{abstract}
\maketitle

\renewcommand{\theitheorem}{\Alph{itheorem}}
\section*{Introduction}
Consider the following setup: $\left(\Omega, \mu\right)$ is a
space with a probability measure $\mu$, and 
$\omega_1, \dots, \omega_n$ is a collection of measurable subsets of
$\Omega$, with $\mu(\omega_i) = p_i.$ We play a game as follows: 
The $j$th step consists of picking a point $x_j \in \Omega$ at random, 
so that after $k$ steps we have the set
$X_k = \{x_1, \dots, x_k\}.$ The game is considered to be over when 
$$
\forall i \leq n, \quad X_k \cap \omega_i \neq X_k.
$$
We consider the duration of our game to be a random variable
$T = T(p_1, \dots, p_n),$ and wish to compute the expection $E(p_1,
\dots, p_n).$
of $T.$ This cannot, in general, be done without knowing the measures
$p_{i_1 i_2 \dots i_k} = \mu(\omega_{i_1} \cap \omega_{i_2} \cap
\cdots \cap \omega_{i_k}),$ and in the sequel we will introduce the 

\medskip\noindent
\textbf{Independence Hypothesis:} 
$$p_{i_1 i_2 \dots i_k} = p_{i_1} \times \dots \times p_{i_k}.$$
Estimates without using the independence hypothesis are shown in 
the companion paper \cite{batch}.

We now assume further that we don't actually know the measures $p_1,
\dots, p_n,$ but know that they themselves are a sample from some (known)
probability distribution $\mathcal{F}$, of necessity supported in $[0,
1].$ We consider $E(p_1, \dots, p_n)$ as our random variable, and we
wish to compute \emph{its} expectation (over the space of all
$n$-element samples from $\mathcal{F},$ and in particular we are
interested in the limiting situation when $n$ is large.

Under the independence assumption, it turns out that we can write
(Lemma \ref{subsuml}): 
\begin{equation}
\label{subsum0} T = \sum_{\mathbf{s}\subseteq \{1, \dots, n\}} 
(-1)^{|s|-1} \left(\frac{1}{1-p_{\mathbf{s}}} - 1\right),
\end{equation}
where if $\mathbf{s} = \{i_1, \dots, i_k\},$ we write 
$p_{\mathbf{s}} = p_{i_1} \times \dots \times p_{i_k}.$
To use equation (\ref{subsum0}) to understand the statistical behavior
of $T$, we must introduce the
\emph{moment zeta function} of the probability distribution
$\mathcal{F},$ defined as follows:
\begin{idefinition}
Let $m_k = \int_0^1 x^k d \mathcal{F}$ be the $k$-th moment of
$\mathcal{F}.$ Then 
$$
\zeta_{\mathcal{F}}(s) = \sum_{k=1}^\infty m_k^s;
$$
\end{idefinition}
The sum in the definition above obviously converges only in some
half-plane $\Re s > s_0$; the function can be analytically continued,
but in the sequel we will be interested in asymptotic results for $s$
a large real number, so this will not use complex variable methods at
all. 

The relevance of this to our questions comes from Lemma
\ref{zetalemma}, which we restate for convenience as:
\begin{ilemma}
\label{izetalemma} Let $\mathcal{F}$ be a probability 
distribution as above, and let $x_1, \dots, x_n$ be independent 
random variables with common distribution $\mathcal{F}$. Then
\begin{equation}
\mathbb{E}\left(\frac{1}{1-x_1 \dots x_n}\right) = 
\zeta_{\mathcal{F}}(n).
\end{equation}
In particular, the expectation is undefined whenever the zeta 
function is undefined. 
\end{ilemma}
Now, we can write (using Lemma \ref{izetalemma}) the following
\emph{formal} identity:
\begin{equation}
\label{formal}
\mathbb{E}(T) = - \sum_{k=1}^n (-1)^k \binom{n}{k} \zeta_{\mathcal{F}}(k).
\end{equation}
The identity is only formal, because $\zeta_{\mathcal{F}}(k)$ is not
necessarily defined for all positive integers $k.$ It \emph{is}
defined for all positive integers $k$ when $\mathcal{F}([1-x, 1]) \sim
x^\alpha,$ for $\alpha > 1$ -- this case is analyzed in Section
\ref{isdef}. If $\alpha = 1$ (we will not deal with the case $\alpha <
1$ in this paper; see \cite{batch}), we write 
\begin{equation*}
T = T_1 - T^\prime,
\end{equation*}
where 
\begin{equation*}
T_1 = \sum_{i=1}^n \frac{1}{1-p_i}.
\end{equation*}
$T_1$ has no expectation, but a $n$ goes to $\infty,$ $T_1/n$ does
converge in distribution to a stable law of exponent $1$ (see
\cite{kr1} and \cite{kr2} for many related results).
The variable $T^\prime$ does possess a finite expectation, given by
\begin{equation}
\label{formal2}
\mathbb{E}(T^\prime) =  \sum_{k=2}^n (-1)^k \binom{n}{k} \zeta_{\mathcal{F}}(k).
\end{equation}
The expressions given by Eq.~(\ref{formal}) and Eq.~(\ref{formal2})
are analyzed in Sections \ref{isdef} and \ref{medalpha}, and the
following Theorems are shown:
\begin{itheorem}[Thm. \ref{mainisdef}]
Let $\mathcal{F}$ be a continuous distribution supported on $[0, 1],$
and let $f$ be the density of $\mathcal{F}.$ Suppose further that 
$$\lim_{x \rightarrow 1} \frac{f(x)}{(1-x)^{\beta}} = c,$$ for $\beta,
c > 0.$ Then, 
\begin{equation*}
\begin{split}
\lim_{n\rightarrow \infty} n^{-\frac{1}{1+\beta}} \left[\sum_{k=1}^n
\binom{n}{k}(-1)^k \zeta_{\mathcal{F}}(k)\right] \\= 
-\int_0^\infty
\frac{1-\exp\left(-c\Gamma(\beta+1)u^{1+\beta}\right)}{u^2} du\\
= - \left(c \Gamma(\beta + 1)\right)^{\frac{1}{\beta+1}}
\Gamma\left(\frac{\beta}{\beta + 1}\right).
\end{split}
\end{equation*}
\end{itheorem}
and 
\begin{itheorem}[Thm. \ref{alpha1asymp}]
\label{ialpha1asymp}
Let $\mathcal{F}$ be a continuous distribution supported on $[0, 1],$
and let $f$ be the density of $\mathcal{F}.$ Suppose further that 
$$\lim_{x \rightarrow 1} \frac{f(x)}{(1-x)} = c > 0.$$ 
Then, 
$$\sum_{k=2}^n
\binom{n}{k}(-1)^k \zeta_{\mathcal{F}}(k) \sim c n \log n.
$$
\end{itheorem}
To get error estimates, we need stronger assumption on the function
$f$ than the weakest possible assumption made in Theorem
\ref{alpha1asymp}. The proof of the below follows by modifying
slightly the proof of Lemma \ref{isest4}:

\begin{itheorem}[Thm. \ref{alpha1asymp2}]
\label{ialpha1asymp2}
Let $\mathcal{F}$ be a continuous distribution supported on $[0, 1],$
and let $f$ be the density of $\mathcal{F}.$ Suppose further that 
$$f(x) \sim c (1-x) + O\left((1-x)^\delta\right),$$ where $\delta > 0.$
Then, 
$$\sum_{k=2}^n
\binom{n}{k}(-1)^k \zeta_{\mathcal{F}}(k) \sim c n \log n + O(n).
$$
\end{itheorem}

Our original probabilistic problem is thus completely resolved, but
the sums given by Eqs. (\ref{formal}) and (\ref{formal2}) are
interesting in and of itself, and, with some more work (Section
\ref{riemannz}), we can
considerably strengthen them as follows for the Riemann zeta function
and its scaling:

\begin{itheorem}[Thm. \ref{riemann}]
\label{iriemann}
$$\sum_{k=2}^n
\binom{n}{k}(-1)^k \zeta(k) \sim  n \log n + (2\gamma-1) n +
O\left(\frac{1}{n}\right),
$$
 where $\zeta$ is the Riemann zeta function
and $\gamma$ is Euler's constant.
\end{itheorem}
and
\begin{itheorem}[Thm. \ref{riemanns}]
\label{iriemanns}
Let $s > 1,$ then
\begin{equation*}
\sum_{k=1}^n
\binom{n}{k}(-1)^k \zeta(s k) \sim \Gamma\left(1-\frac{1}{s}\right) n^{\frac{1}{s}}.
\end{equation*}
\end{itheorem}

It should be remarked that using the methods of section \ref{riemannz}
higher order terms in the asymptotics can be obtained, if
desired, but they seem to be of more limited interest.

\section{A formula for the winning time $T$}
An application of the 
inclusion-exclusion principle gives us:

\thm{lemma}{\label{latmost} The probability $l_k$ that we have 
won after $k$ steps is given by 
$$
l_k=\prod_{i=1}^n(1-p_i^k).
$$
}

Note that the probability $s_k$ of winning the game \emph{on the 
$k$-th step} is given by $s_k = l_k - l_{k-1}= (1-l_{k-1}) - 
(1-l_k)$. Since the expected number of steps $T$ 
is given by
$$T = \sum_{k=1}^\infty k s_k,$$
we immediately have  $$T = \sum_{k=1}^\infty (1-l_k)$$
\begin{lemma}
\label{letime} The expected time $T$ of learning the 
concept $R_0$ is given by
\begin{equation}
\label{letimeeq}
T = \sum_{k=1}^\infty \left(1-\prod_{i=1}^n 
\left(1-p_i^k\right)\right).
\end{equation}
\end{lemma}
Since the sum above is absolutely convergent, we can expand the 
products and interchange the order of summation to get the 
following formula for $T$:

\medskip\noindent
\textbf{Notation.}
Below, we identify subsets of $\{1, \dots, n\}$ with 
 multindexes (in the obvious way), and if $s = \{i_1, \dots, i_l\},$ then
$$p_s \stackrel{\mbox{def}}= p_{i_1} \cdots p_{i_l}.$$

\begin{lemma}
\label{subsuml}
The expression Eq. (\ref{letimeeq}) can be rewritten as:
\begin{equation}
\label{subsum} T = \sum_{s\subseteq \{1, \dots, n\}} 
(-1)^{|s|-1} \left(\frac{1}{1-p_s} - 1\right),
\end{equation}
\end{lemma}

\begin{proof}
With notation as above,
\begin{equation*}
\prod_{i=1}^m \left(1-p_i^k\right) = 
\sum_{s \subseteq \{1, \dots, n\}} (-1)^{|s|} p_s^k,
\end{equation*}
so
\begin{equation*}
\begin{split}
T &= \sum_{k=1}^\infty \left(1 - \prod_{i=1}^n 
\left(1-p_i^k\right)\right)\\ 
&= \sum_{k=1}^\infty \left(1-\sum_{s \subseteq \{1, \dots, n\}}
(-1)^{|s|} p_s^k\right)\\
&= \sum_{s\subseteq \{1, \dots, n\}} (-1)^{|s|-1} 
\sum_{k=1}^\infty p_s^k \\
&= \sum_{s\subseteq \{1, \dots, n\}} 
(-1)^{|s|-1} \left(\frac{1}{1-p_s} - 1\right),
\end{split}
\end{equation*}
where the change in the order of summation is permissible since all
sums converge absolutely.
\end{proof}
Formula (\ref{subsum}) is useful in and of itself, but we now 
use it to analyse the statistical properties of the time of 
success $T$ under our distribution and independence assumptions. 
For this we shall need to study the \emph{moment zeta function} of a
probability distribution, introduced below.
\section{Moment zeta function}
\begin{definition}
\label{zdef} Let $\mathcal{F}$ be a probability 
distribution on a (possibly infinite) interval $I$, and let
$m_k(\mathcal{F}) =  \int_I x^k\mathcal{F}(d x)$ be the $k$-th moment
of  $\mathcal{F}$. Then the \emph{moment zeta function of  
$\mathcal{F}$} is defined to be $$\zeta_{\mathcal{F}}(s) =
\sum_{k=1}^\infty m_k^s(\mathcal{F}),$$ whenever the sum is defined. 
\end{definition}
The definition is, in a way, motivated by the following:

\begin{lemma}
\label{zetalemma} Let $\mathcal{F}$ be a probability 
distribution as above, and let $x_1, \dots, x_n$ be independent 
random variables with common distribution $\mathcal{F}$. Then
\begin{equation}
\mathbb{E}\left(\frac{1}{1-x_1 \dots x_n}\right) = 
\zeta_{\mathcal{F}}(n).
\end{equation}
In particular, the expectation is undefined whenever the zeta 
function is undefined. 
\end{lemma}
\begin{proof}
Expand the fraction in a geometric series and apply Fubini's 
theorem.
\end{proof}
\begin{example}
For $\mathcal{F}$ the uniform distribution on 
$[0, 1]$, $\zeta_{\mathcal{F}}$ is the familiar Riemann zeta 
function. 
\end{example}

Our first observation is that for distributions supported in $[0, 1]$,
the asymptotics of the  moments are determined by the local properties of the 
distribution at $x=1.$ To show this, first recall that 
the \emph{Mellin transform} of $f$ is defined to be
$$\mathcal{M}(f)(s) = \int_0^1 f(x) x^{s-1} d x.$$ 
Mellin transform is closely related to the Laplace
transform. Making the  
substitution $x = \exp(-u)$, we see that $$\mathcal{M}(f) = 
\int_0^\infty f(\exp(-u)) \exp(-s u) d u,$$ so the Mellin 
transform of $f$ is equal to the Laplace transform of $f \circ 
\exp,$ where $\circ$ denotes functional composition.

The following observation is both obvious and well-known:
\begin{lemma}
$m_k(\mathcal{F}) = \mathcal{M}(f)(k+1).$
\end{lemma}
It follows that computing the asymptotic behavior of the $k$-th moment
of $\mathcal{F}$ as a function of $k$ reduces to calculating the large
$s$ asymptotics of the Mellin transform, which is tantamount to
computing the asymptotics of the Laplace transform of $f \circ \exp.$

\begin{theorem}
\label{momasymp}
Let $\mathcal{F}$ be a continuous distribution supported in $[0, 1],$ 
let $f$ be the density of the distribution $\mathcal{F}$, and
suppose that $f(1-x) = c x^\beta + O(x^{\beta + \delta}),$ for some
$\delta > 0.$ Then the $k$-th moment of $\mathcal{F}$ is asymptotic to 
$C k^{-(1+\beta)},$ for $C = c \Gamma(\beta+1).$ 
\end{theorem}

\begin{proof} 
The asymptotics of the Laplace transform are easily 
computed by Laplace's method, and in the case we are interested 
in, Watson's lemma (see, eg, \cite{benorsz}) tells us that if 
$f(x) \asymp c (1-x)^\beta$, then $\mathcal{M}(f)(s) \asymp c 
\Gamma(\beta+1) x^{-(\beta + 1)}.$ 
\end{proof}
\begin{corollary}
Under the assumptions of Theorem \ref{momasymp}, 
$\zeta_{\mathcal{F}}(s)$ is defined for $s 
>1/(1+\beta)$.
\end{corollary}
We will need another observation:
\begin{lemma}
\label{decmom}
For $\mathcal{F}$ supported in $[0, 1]$, $m_k(\mathcal{F})$ is
monotonically decreasing as a function of $k.$
\end{lemma}
\begin{proof}
Immediate.
\end{proof}
Below we shall analyze three cases.
In the sequel, we set $\alpha = \beta + 1$.
\section{$\alpha > 1$}
\label{isdef} 
In this case, we use our assumptions to rewrite Eq. 
(\ref{subsum}) as 
\begin{equation}
\label{subsum2} 
\mathbb{E}(T) = - \sum_{k=1}^n \binom{n}{k}(-1)^k \zeta_{\mathcal{F}}(k).
\end{equation}
This, in turn, can be rewritten (by expanding the definition of 
zeta) as
\begin{equation}
\label{subsum3} \mathbb{E}(T) = - \sum_{j=1}^\infty 
\left[\left(1-m_j(\mathcal{F})\right)^n-1\right] = 
\sum_{j=1}^\infty \left[1- \left(1-m_j(\mathcal{F})\right)^n\right]
\end{equation}
Since the terms in the sum is monotonically decreasing (as a function
of $j$) by Lemma \ref{decmom}, the sum in 
Eq. (\ref{subsum3}) can be approximated by an integral of 
\emph{any} monotonic interpolation $m$ of the sequence 
$m_j(\mathcal{F})$ -- we will interpolate by $m(x) = 
\mathcal{M}(f)(x+1)$). The error of such an approximation is bounded
by the first term,  which, in turn, isbounded in absolute value by $2$, to get 
\begin{equation}
\label{approx1} T = - \int_1^\infty \left[(1-m(x))^n -1\right] d 
x + O(1),
\end{equation}
where the error term is bounded above by $2$.
We shall write
\begin{equation*}
T_0 = - \int_1^\infty \left[(1-m(x))^n -1\right]  d x.
\end{equation*}
Now, let us assume that 
\begin{equation}
\label{asympass}
\lim_{x\rightarrow \infty} x^{\alpha} m(x)=L,
\end{equation}
 for some 
$\alpha > 1.$ We substitute $x = n^{1/\alpha}/u$, to get
\begin{equation*}
 T_0 =  n^{\frac{1}{\alpha}}\int_0^{n^{\frac{1}{\alpha}}}
\frac{\left[1-\left(1-m(n^{1/\alpha}/u)\right)^n \right]}{u^2} d u = n^{\frac{1}{\alpha}}\left[I_1(n) +
I_2(n)\right],
\end{equation*}
where
\begin{equation*}
I_1(n) = \int_0^{n^{\frac{1}{3 \alpha}}}
\frac{\left[1-(1-m(n^{1/\alpha}/u)^n \right]}{u^2} d u,
\end{equation*}
and 
\begin{equation*}
I_2(n) = \int_{n^{\frac{1}{3 \alpha}}}^{n^{\frac{1}{\alpha}}}
\frac{\left[1-(1-m(n^{1/\alpha}/u)^n \right]}{u^2} d u.
\end{equation*}
We will need the following:
\begin{lemma}
\label{explem}
Let $f_n(x) = (1-x/n)^n,$ and let $0 \leq z < 1/2.$ 
$$f_n(x) = \exp(-x)\left[1-\frac{x^2}{2 n} + O\left(\frac{x^3}{n^2}\right)\right].$$
\end{lemma}
\begin{proof}
Note that 
$$\log f_n(x) = n \log(1-x/n) = -x - \sum_{k=2}^\infty \frac{x^k}{kn^{k-1}}.$$
The assertion of the lemma follows by exponentiating the two sides of
the above equation.  
\end{proof}
Now:
\begin{lemma}
\label{assi2}
\begin{equation*}
\lim_{n \rightarrow \infty} n^{\frac{1}{\alpha}} I_2(n) = 0.
\end{equation*}
\end{lemma}
\begin{proof}
The integrand of $I_2(n)$ is monotonically decreasing, and so 
$$I_2(n) \leq n^{-\frac{2}{3 \alpha}}
\left[1-\left(1-m\left(n^{\frac{2}{3 \alpha}}\right)\right)^n\right].$$
By our assumption Eq. (\ref{asympass}) and by Lemma \ref{explem} we
see that the right hand side goes to zero (exponentially fast).
\end{proof}
\begin{lemma}
\label{assi1}
$$\lim_{n\rightarrow \infty} I_1(n) = \int_0^\infty
\frac{1-\exp\left(-Lu^{\alpha}\right)}{u^2} d u.$$
\end{lemma}
\begin{proof}
Immediate from Eq.~(\ref{asympass}) and Lemma \ref{explem}. Note that
the integral converges when $\alpha$ is greater than $1.$
\end{proof}
\begin{remark}
\label{euler0}
$$\int_0^\infty
\frac{1-\exp\left(-Lu^{\alpha}\right)}{u^2} d u = L^{\frac{1}{\alpha}}
\Gamma\left(\frac{\alpha - 1}{\alpha}\right).$$
\end{remark}
\begin{proof}
\begin{equation*}
\int_0^\infty
\frac{1-\exp\left(-Lu^{\alpha}\right)}{u^2} d u =
\lim_{\epsilon \rightarrow 0} \left[\frac{1}{\epsilon} - 
\int_\epsilon^\infty \frac{\exp\left(-L u^{\alpha}\right)}{u^2} d u\right].
\end{equation*}
To prove the remark we need to analyze the behavior of the integral
above as $\epsilon \rightarrow 0.$ First, we change variables:
$v = L u^\alpha.$ Then,
\begin{equation*}
\int_\epsilon^\infty \frac{\exp\left(-L u^{\alpha}\right)}{u^2} d u = 
\frac{L^{1/\alpha}}{\alpha}\int_{L \epsilon^\alpha}^\infty
\exp(-v)v^{-(1+1/\alpha)}
d v.
\end{equation*}
Integrating by parts, get
\begin{equation*}
\int_{L \epsilon^\alpha}^\infty \exp(-v)v^{-(1+1/\alpha)} d v =
-\alpha \left. \exp(-v) v^{1/\alpha} \right|_{L \epsilon^\alpha}^\infty - 
\alpha \int_{L \epsilon^\alpha} \exp(-v) v^{-1/\alpha}.
\end{equation*}
Since $1/\alpha < 1,$ $\int_0^\infty \exp(-v) v^{-1/\alpha} d v = 
\Gamma(1-1/\alpha),$ from which the assertion of the remark follows.
\end{proof}
We summarize as follows:
\begin{theorem}
\label{mainisdef}
Let $\mathcal{F}$ be a continuous distribution supported on $[0, 1],$
and let $f$ be the density of $\mathcal{F}.$ Suppose further that 
$$\lim_{x \rightarrow 1} \frac{f(x)}{(1-x)^{\beta}} = c,$$ for $\beta,
c > 0.$ Then, 
\begin{equation*}
\begin{split}
\lim_{n\rightarrow \infty} n^{-\frac{1}{1+\beta}} \left[\sum_{k=1}^n
\binom{n}{k}(-1)^k \zeta_{\mathcal{F}}(k)\right] \\= 
-\int_0^\infty
\frac{1-\exp\left(-c\Gamma(\beta+1)u^{1+\beta}\right)}{u^2} du\\
= - \left(c \Gamma(\beta + 1)\right)^{\frac{1}{\beta+1}}
\Gamma\left(\frac{\beta}{\beta + 1}\right).
\end{split}
\end{equation*}
\end{theorem}
\begin{proof}
Follows from Lemmas \ref{assi1} and \ref{assi2} together with
theorem \ref{momasymp} and Remark \ref{euler0}.
\end{proof}

\section{$\alpha = 1$}
\label{medalpha} In this case, 
\begin{equation}
\label{asest02}
f(x) = L + o(1)
\end{equation} as $x$ 
approaches $1,$ and so Theorem \ref{momasymp} tells us that 
\begin{equation}
\label{asest2}
\lim_{j \rightarrow \infty} j m_j(\mathcal{F}) = L.
\end{equation}
It is not hard to see that 
$\zeta_{\mathcal{F}}(n)$ is defined for $n \geq 2$. We break up 
the expression in Eq. (\ref{subsum}) as 
\begin{equation}
\label{subsumm} T = \sum_{j=1}^n {\frac{1}{1-p_j} - 1} + 
\sum_{s\subseteq \{1, \dots, n\}, \quad |s| > 1} 
 (-1)^{|s|-1} 
\left(\frac{1}{1-p_s} - 1\right).
\end{equation}
Let 
\begin{gather*} T_1 = \sum_{j=1}^n {\frac{1}{1-p_j} - 1},\\
 T_2 = \sum_{s\subseteq \{1, \dots, n\}, \quad |s| > 1} 
 (-1)^{|s|-1} 
\left(\frac{1}{1-p_s} - 1\right).
\end{gather*}
 The first sum $T_1$ has 
no expectation, however $T_1/n$  does have have a stable 
distribution centered on $c \log n + c_2$. We will keep this in 
mind, but now let us look at the second sum  $T_2$. It can be 
rewritten as 
\begin{equation}
\label{subsumm2} T_2(n) = - \sum_{j=1}^\infty 
\left[\left(1-m_j(\mathcal{F})\right)^n-1 + n
m_j(\mathcal{F})\right]. 
\end{equation}
\begin{lemma}
\label{monoton}
The quantity $y_j = \left(1-m_j(\mathcal{F})\right)^n-1 + n
m_j(\mathcal{F})$ is a monotonic function of $j.$
\end{lemma}
\begin{proof}
We know that $m_j(\mathcal{F})$ is a monotonically decreasing
positive function of $j,$ and that $m_0(\mathcal{F}) = 1.$ It is
sufficient to show that the function  $g_n(x) = (1-x)^n + n x$ is
monotonic for $x \in (0, 1].$ We compute
$$\frac{d g_n(x)}{d x} = n \left(1-(1-x)^{n-1}\right) > 0,$$
for $x \in (0, 1).$
\end{proof}
Lemma \ref{monoton} allows us to use the same method as in section
\ref{isdef} under the assumption  
that the $k$-th moment is asymptotic to $k^\alpha$ (this time for 
$\alpha \leq 1$). Since the term $y_j$ is bounded above by a constant times $n$, we can write
\begin{equation}
\label{zetabreak}
T_2(n) = S_2(n) + O(n),
\end{equation}
where
\begin{equation}
S_2(n) = n \int_0^n \frac{\left[1-n m(n/u) - (1-m(n/u)^n 
\right]}{u^2} d u.
\end{equation}
\begin{remark}
\label{zetab2}
The error term in Eq. (\ref{zetabreak}) above can be improved in the
case where $\mathcal{F}$ is the uniform distribution on $[0, 1],$
in which case $m_j = 1/j.$ In that case $T_2(n) = S_2(n) - \gamma n +
O(1),$ where $\gamma$ is Euler's constant.
\end{remark}
\begin{proof}
In this case, we write 
\begin{equation*}
\begin{split}
T_2(n) = \lim_{k\rightarrow\infty} 
- \sum_{j=1}^k
\left[\left(1-m_j(\mathcal{F})\right)^n-1 + n
m_j(\mathcal{F})\right]\\ = - \sum_{j=1}^k
\left[\left(1-m_j(\mathcal{F})\right)^n-1\right]
-n \sum_{j=1}^k m_j.
\end{split}
\end{equation*}
The terms in the first sum are decreasing, so the
first sum can be approximated by an integral with total error $O(1).$
As for the second sum, since $m_j = 1/j,$ it is well-known (eg,
Euler-Maclaurin summation) that 
$$\sum_{j=1}^k \frac{1}{j} = \int_1^k \frac{dx}{x} + \gamma +
O(\frac{1}{k}),$$
from which the assertion of the remark follows.
\end{proof}
To understand the asymptotic behavior of $S_2(n)$ we write
$$S_2(n) = n \left[I_1(n) + I_2(n) + I_3(n) + I_4(n)\right],$$
where
\begin{gather}
I_1(n) = \int_0^1 \frac{\left[1-n m(n/u) -
\left(1-m\left(\frac{n}{u}\right)\right)^n\right]}{u^2} d u,\\
I_2(n) = \int_1^{n^{\frac13}} \frac{\left[1-
\left(1-m\left(\frac{n}{u}\right)\right)^n\right]}{u^2} d u,\\
I_3(n) = \int_{n^{\frac13}}^n \frac{\left[1-
\left(1-m\left(\frac{n}{u}\right)\right)^n\right]}{u^2} d u,\\
I_4(n) = -n \int_1^n \frac{m\left(n/u\right)}{u^2}
d u.
\end{gather}
\begin{lemma}
\label{intest1}
$$\lim_{n\rightarrow \infty} I_1(n) =
 \int_0^1 \frac{1-\exp\left(-L u\right) - L u}{u^2}d u.$$
\end{lemma}
\begin{proof}
Immediate from the estimate Eq.~(\ref{asest2}) and Lemma \ref{explem}.
\end{proof}
\begin{lemma}
\label{intest2}
$$\lim_{n\rightarrow \infty} I_2(n) =
 \int_1^\infty \frac{1-\exp\left(-L u\right)}{u^2} d u.$$
\end{lemma}
\begin{proof}
Again, immediate from Eq.~(\ref{asest2}) and Lemma \ref{explem}.
\end{proof}
\begin{remark}
\label{euler}
$$\int_0^1 \frac{1-\exp\left(-L u\right) - L u}{u^2} +  \int_1^\infty
\frac{1-\exp\left(-L u\right)}{u^2} = L(1-\gamma - \log L),$$
where $\gamma$ is Euler's constant.
\end{remark}
\begin{proof}
\begin{equation}
\begin{split}
\label{firsthalf}
\int_0^1 \frac{1-\exp\left(-L u\right) - L u}{u^2}d u +  \int_1^\infty
\frac{1-\exp\left(-L u\right)}{u^2} d u =\\
\lim_{\epsilon \rightarrow 0} \left\{
- \int_\epsilon^\infty \left[\frac{\exp(-L u)}{u^2} + \frac{1}{u^2}\right]
d u -  L \int_\epsilon^1 \frac{du}{u}\right\} = \\
\lim_{\epsilon \rightarrow 0} \left\{\frac{1}{\epsilon} +  L \log \epsilon -
\int_{\epsilon}^\infty \frac{\exp(-L u)}{u^2} d u\right\}.
\end{split}
\end{equation}
To evaluate the last limit, we need to compute the expansion as
$\epsilon \rightarrow 0$ of the last integral. Changing variables
$v = L u,$ we get
\begin{equation*}
\begin{split}
\int_\epsilon^\infty
\frac{\exp(-Lu)}{u^2} d u =
L \int_{L\epsilon}^\infty \frac{\exp(-v)}{v^2} d v = \\
L \left[ -\left. 
\frac{\exp(-v)}{v} \right|_{L \epsilon}^\infty
 - \int_{L \epsilon}^\infty \frac{\exp(-v)}{v} d v \right] = \\
L \left[\frac{\exp(-v)}{L \epsilon} - \left.\exp(-v) \log(v)
\right|_{L\epsilon}^\infty 
- \int_{L \epsilon}^\infty \exp(-v) \log(v) dv\right]=\\
\frac{\exp(- L \epsilon)}{\epsilon} + 
L \exp(-L\epsilon)\log(\epsilon) +\\ 
L \log L \exp(- L \epsilon) - 
L \int_{L\epsilon}^\infty 
\exp(-v) \log(v) d v.
\end{split}
\end{equation*}
Substituting into Eq.~(\ref{firsthalf}), we get
\begin{equation*}
\begin{split}
\int_0^1 \frac{1-\exp\left(-L u\right) - L u}{u^2}d u +  \int_1^\infty
\frac{1-\exp\left(-L u\right)}{u^2} d u =\\
\lim_{\epsilon \rightarrow 0} 
\left\{\frac{1-\exp(-L \epsilon)}{\epsilon} 
+ L(1-\exp(-L \epsilon))
\log \epsilon \right.\\ \left. - L \log L \exp(- L \epsilon) + \int_{L \epsilon}^\infty
\exp(-v) \log v d v\right\}\\
= L\left(1 - \log L + \int_{L \epsilon}^\infty \exp(-v) \log v dv\right).
\end{split}
\end{equation*}
Since $\int_0^\infty \log(x) \exp(-x) d x = - \gamma,$  the result follows.
\end{proof}
\begin{lemma}
\label{isest3}
$$\lim_{n \rightarrow \infty} n I_3(n) = 0.$$
\end{lemma}
\begin{proof}
See the proof of Lemma \ref{assi2}.
\end{proof}
\begin{lemma}
\label{isest4}
$$\lim_{n \rightarrow \infty}
-\frac{I_4(n)}{\log n} = L.$$
\end{lemma}
\begin{proof}
We shall show that the limit in question lies between $(1-\epsilon)L$
and $(1+\epsilon)L,$ for any $\epsilon > 0,$ from which the conclusion
of the lemma obviously follows. To do that, pick $C,$ such that 
$$1-\epsilon/4 \leq x  m(x) \leq 1+\epsilon/4$$
for $x > C.$
Now, write
$$\int_1^n \frac{m(n/u)}{u^2} d u = J_1(n) + J_2(n),$$
where
\begin{gather}
J_1(n) = \int_1^{\frac{n}{C}} \frac{m(n/u)}{u^2} d u\\
J_2(n) = \int_{\frac{n}{C}}^n \frac{m(n/u)}{u^2} d u.
\end{gather}
Observe that 
$$0 < J_2(n) = \frac{1}{n}\int_1^C m(x) d x \leq \frac{C-1}{n},$$
while
$$\frac{1-\epsilon/4}{n} \int_1^{\frac{n}{C}} \frac{d u}{u}
\leq J_1(n) \leq \frac{1+\epsilon/4}{n} \int_1^{\frac{n}{C}}
\frac{d u }{u},
$$
so
$$\frac{(1-\epsilon/4)}{n} (\log n - \log C) \leq J_1(n) \leq 
\frac{(1+\epsilon/4)}{n} (\log n - \log C).$$
If we now pick $N = C^{4/\epsilon},$ it is clear that for $n > N,$
$$(1-\epsilon/2) \log n \leq J_1(n) \leq (1+\epsilon/2) \log n,$$
while $J_2$ is bounded above in absolute value by $C^{1-4/\epsilon}.$
\end{proof}
The above lemmas can be summarized in the following
\begin{theorem}
\label{alpha1asymp}
Let $\mathcal{F}$ be a continuous distribution supported on $[0, 1],$
and let $f$ be the density of $\mathcal{F}.$ Suppose further that 
$$\lim_{x \rightarrow 1} \frac{f(x)}{(1-x)} = c > 0.$$ 
Then, 
$$\sum_{k=2}^n
\binom{n}{k}(-1)^k \zeta_{\mathcal{F}}(k) \sim c n \log n.
$$
\end{theorem}
To get error estimates, we need stronger assumption on the function
$f$ than the weakest possible assumption made in Theorem
\ref{alpha1asymp}. The proof of the below follows by modifying
slightly the proof of Lemma \ref{isest4}:

\begin{theorem}
\label{alpha1asymp2}
Let $\mathcal{F}$ be a continuous distribution supported on $[0, 1],$
and let $f$ be the density of $\mathcal{F}.$ Suppose further that 
$$f(x) \sim c (1-x) + O\left((1-x)^\delta\right),$$ where $\delta > 0.$
Then, 
$$\sum_{k=2}^n
\binom{n}{k}(-1)^k \zeta_{\mathcal{F}}(k) \sim c n \log n + O(n).
$$
\end{theorem}
\section{Riemann zeta function}
\label{riemannz}
The proof of the key Lemma \ref{isest4} is trivial in the case where
$f(x) = 1,$ and so $\zeta_{\mathcal F}$ is the Riemann zeta
function. In that case, however, we get the following much stronger
result:
\begin{theorem}
\label{riemann}
$$\sum_{k=2}^n
\binom{n}{k}(-1)^k \zeta(k) \sim  n \log n + (2\gamma-1) n +
O\left(\frac{1}{n}\right),
$$
 where $\zeta$ is the Riemann zeta function
and $\gamma$ is Euler's constant.
\end{theorem}
It should also be noted that the results of Section \ref{isdef} immediately
imply the following:
\begin{theorem}
\label{riemanns}
Let $s > 1,$ then
\begin{equation*}
\sum_{k=1}^n
\binom{n}{k}(-1)^k \zeta(s k) \sim \Gamma\left(1-\frac{1}{s}\right) n^{\frac{1}{s}}.
\end{equation*}
\end{theorem}
To prove Theorem \ref{riemann} we need to sharpen some of the
estimates of the preceding section. First:
\begin{lemma}
\label{rintest}
Let the notation be as in the preceding section.
When $m(x) = \frac{1}{x},$
\begin{gather}
I_1(n) =  \int_0^1 \frac{1-\exp\left(- u\right) -  u}{u^2} d u +
\frac{1}{2 n} \int_0^1 \exp(- u) d u + O\left(\frac{1}{n^2}\right),\\
I_2(n) = \int_1^\infty \frac{1-\exp\left(- u\right)}{u^2} d u +
\frac{1}{2 n} \int_1^\infty \exp(- u) d u.
\end{gather}
\end{lemma}

\begin{proof}
Imediate from the expansion in Lemma \ref{explem}.
\end{proof}
We can also sharpen the statement of Lemma \ref{isest3}:
\begin{lemma}
\label{rintest3}
\begin{equation*}
\lim_{n\rightarrow \infty} n^k I_3(n) = 0,
\end{equation*}
for any $k.$
\end{lemma}
\begin{proof}
This statement holds in general, and no change in argument is
necessary.
\end{proof}
In the case where $m(x) = 1/x,$ Lemma \ref{isest4} is immediate, and
has no error term:
\begin{lemma}
\label{rintest4}
\begin{equation*}
I_4(n) = - \log(n).
\end{equation*}
\end{lemma}
\begin{proof} Immediate. \end{proof}
We now have:
\begin{theorem}
\label{bigothm}
$$\sum_{k=2}^n
\binom{n}{k}(-1)^k \zeta(k) \sim  n \log n + (2\gamma-1) n + O\left(1\right),
$$
\end{theorem}
\begin{proof}
Lemmas \ref{rintest}, \ref{rintest3} and \ref{rintest4}, combined with
Remark \ref{zetab2}.
\end{proof}

To improve the error term from that in Theorem \ref{bigothm}, it is
necessary to sharpen the estimate in Remark \ref{zetab2} to:
\begin{theorem} 
\label{sumint}
With the notation of Remark \ref{zetab2}, 
\begin{equation*}
T_2(n) = S_2(n) - \gamma n - \frac{1}{2} + O\left(\frac{1}{n}\right).
\end{equation*}
\end{theorem}
\begin{proof}
The Theorem will follow immediately from Lemma \ref{si1} and the
results of Section \ref{sumintsec}.
\end{proof}
\begin{lemma}
\label{si1}
\begin{equation*}
\lim_{N \rightarrow \infty}\sum_{j=1}^N \frac{1}{j}- \log(n) = \gamma.
\end{equation*}
\end{lemma}
\begin{proof}
Well-known.
\end{proof}
\subsection{A sum and an integral}
\label{sumintsec}
Let
\begin{gather*}
S_n(N) = \sum_{j=1}^N \left(1 - \frac{1}{j}\right)^n,\\
I_n(N) = \int_1^N \left(1 - \frac{1}{x}\right)^n d x,\\
D_n(N) = S_n(N) - I_n(N),\\
D_n = \lim_{N \rightarrow \infty} D_n(N).
\end{gather*}
In this section we shall prove the following result:
\begin{theorem}
\label{sidif}
\begin{equation*}
D_n = \frac{1}{2} + o\left(\frac{1}{n}\right).
\end{equation*}
\end{theorem}
We will need the following preliminary results:
\begin{lemma}
\label{eumac}
Let $f$ be a $C^1$ function defined on $[0, \infty).$ Then
\begin{equation*}
\sum_{k=0}^N f(k) = \frac{1}{2}\left[f(0) + f(N)\right]
+ \int_0^Nf(t) d t + \int_0^N \left(\left\{t\right\}- \frac{1}{2}\right)
f^\prime(t) d t.
\end{equation*}
\end{lemma}
\begin{proof}
Integration by parts -- see Exercises for Section 6.7 of \cite{benorsz}.
\end{proof}
\begin{lemma}
\label{genlem}
Let $f$ be a $C^2$ function defined on $[0, \infty),$ such that $f^{\prime\prime}$ is bounded, and $f^{\prime\prime}(x) = O(1/x^2).$
Then 
\begin{equation*}
\left|\sum_{k=0}^\infty \int_{\frac{k}{n}}^{\frac{k+1}{n}}
\left(x - \frac{k+\frac{1}{2}}{n}\right) f(x) d x
- \frac{1}{4 n^3} \sum_{k=0}^\infty
f^\prime\left(\frac{k+\frac{1}{2}}{n}\right)\right| =
O\left(\frac{1}{n^4}\right).
\end{equation*}
\end{lemma}
\begin{proof}
On the interval $\left[k/n, (k+1)/n\right]$ we can write 
\begin{equation}
\label{auxgenlem}
f(x) = f\left(\frac{k+\frac{1}{2}}{n}\right) +
f^\prime\left(\frac{k+\frac{1}{2}}{n}\right)
\left(x -\frac{k+\frac{1}{2}}{n}\right) + R_2(x),
\end{equation}
where, by Taylor's theorem, $|R_2(x)| \leq x^2 \max_{x \in [k/n, (k+1)/n]}
f^{\prime\prime}(x).$ The assertion of the lemma then follows by
integration of Eq. (\ref{auxgenlem}).
\end{proof}
\begin{lemma}
\label{genlem2}
Under the assumptions of Lemma \ref{genlem}, together with the
assumtption that $f$ and all of its derivatives vanish at $0$
\begin{equation*}
\left| \sum_{k=0}^\infty f^\prime\left(\frac{k+\frac{1}{2}}{n}\right)\right|
 = O\left(\frac{1}{n}\right).
\end{equation*}
\end{lemma}
\begin{proof}
Let $g(y) = f^\prime((x+ 1/2)/n).$
Then:
\begin{equation*}
\begin{split}
\sum_{k=0}^\infty f^\prime\left(\frac{k+\frac{1}{2}}{n}\right) 
&= \sum_{k=0}^\infty g(k)\\
&= \frac{1}{2}g(0) + \int_0^\infty g(x) d x + \int_0^\infty
\left(\{x\} - \frac12\right)g^\prime(x) d x\\
&= \frac{1}{2}f^\prime\left(\frac{1}{2 n}\right) + \int_0^\infty
f^\prime\left(\frac{x+\frac{1}{2}}{n}\right) d x +
O\left(\frac{1}{n}\right)\\
&= n \int_{\frac{1}{2 n}}^\infty f^\prime(x) d x +
O\left(\frac{1}{n}\right)\\
&= O\left(\frac{1}{n}\right).
\end{split}
\end{equation*}
\end{proof}

Now we proceed to the proof of Theorem \ref{sidif} \emph{per se}.
First:
\begin{lemma}
\label{dnform}
\begin{equation*}
D_n = \frac{1}{2} + n 
\int_1^\infty \frac{\left(\{x\} -
\frac{1}{2}\right)\left(1-\frac{1}{x}\right)^{n-1}}{x^2} d x.
\end{equation*}
\end{lemma}
\begin{proof}
Immediate corollary of Lemma \ref{eumac}.
\end{proof}
\begin{proof}[Proof of Theorem \ref{sidif}]
By Lemma \ref{dnform} it remains to analyse the asymptotic behavior of
\begin{equation*}
J_n = (n+1) \int_1^\infty \frac{\left(\{x\} -
\frac{1}{2}\right)\left(1-\frac{1}{x}\right)^n}{x^2} d x.
\end{equation*}
(the expression occuring in Lemma \ref{dnform} is actually $J_{n-1}$,
we have changed the variable for notational convenience).
First, we make the substitution $x = n y$, to get
\begin{equation*}
J_n = \frac{(n+1)}{n} \underbrace{\int_{\frac{1}{n}}^\infty
\frac{\left(\{n y\} -\frac{1}{2}\right)\left(1-\frac{1}{n
y}\right)^n}{y^2} d x}_{K_n}, 
\end{equation*}
where clearly $J_n \sim K_n.$
We now write 
$$K_n = \left[\underbrace{\int_{\frac{1}{n}}^{n^{-\frac{1}{3}}}}_{K^\prime_n} + 
\underbrace{\int_{n^{-\frac{1}{3}}}^\infty}_{K^{\prime\prime}_n}
\right]
\frac{\left(\{n y\} -\frac{1}{2}\right)\left(1-\frac{1}{n
y}\right)^n}{y^2} d x.
$$
The integrand of $K^\prime_n$ is bounded above by 
\begin{equation*}
\left(1-n^{-\frac{2}{3}}\right)^n,
\end{equation*}
while the interval of integration is polynomial in length, which
implies that $K_n^\prime$ decreases faster than any power of $n,$ and
so can be ignored for our purposes. On the other hand, Lemma
\ref{explem} implies that 
\begin{equation*}
\begin{split}
K^{\prime\prime}_n &\sim \int_0^\infty \left(\{n y\} -
\frac{1}{2}\right)
\frac{\exp\left(-\frac{1}{y}\right)}{y^2} d y\\
&= \sum_{k=0}^\infty \int_{\frac{k}{n}}^{\frac{k+1}{n}}
\left[ n y - \frac{1}{2} - k\right]
\frac{\exp\left(-\frac{1}{y}\right)}{y^2} d y \\
&= n \sum_{k=0}^\infty \int_{\frac{k}{n}}^{\frac{k+1}{n}}
\left[  y - \frac{k+\frac{1}{2}}{n}\right]
\frac{\exp\left(-\frac{1}{y}\right)}{y^2} d y.
\end{split}
\end{equation*}
We can now apply Lemmas \ref{genlem} and \ref{genlem2} with 
$$f(x) = \frac{\exp\left(-\frac{1}{y}\right)}{y^2};$$ it is easy to
check that $f(x)$ satisfies the assumptions. Theorem \ref{sidif}
follows.
\end{proof}

\end{document}